\documentclass{amsart}
\usepackage{amsmath}
\usepackage{graphicx}
\usepackage{latexsym}

\newtheorem{thm}{Theorem}
\newtheorem{lem}[thm]{Lemma}
\newtheorem{cor}[thm]{Corollary}
\newtheorem{prop}[thm]{Proposition}

\theoremstyle{definition}
\newtheorem{defn}[thm]{Definition}
\newtheorem{rk}[thm]{Remark}

\newcommand{\CPb}{\overline{\mathbb{CP}}{}^{2}}
\newcommand{\CP}{{\mathbb{CP}}{}^{2}}

\newcommand{\Z}{\mathbb{Z}}

\def \x {\times}

\begin{document}

\title[Classification of broken Lefschetz fibrations with small fiber genera]
{Classification of broken Lefschetz fibrations with small fiber genera}
\vspace{0.2in} 

\author[R. \.{I}. Baykur]{R. \.{I}nan\c{c} Baykur}
\address{Brandeis University, Department of Mathematics} 
\email{baykur@brandeis.edu }

\author[S. Kamada]{Seiichi Kamada}
\address{Hiroshima University, Department of Mathematics} 
\email{kamada@math.sci.hiroshima-u.ac.jp}

\begin{abstract}
In this article, we generalize the classification of genus one Lefschetz fibrations to genus one simplified broken Lefschetz fibrations, which have fibers of genera one and zero. We classify genus one Lefschetz fibrations over the $2$-disk with certain non-trivial global monodromies using chart descriptions, and identify the $4$-manifolds admitting genus one simplified broken Lefschetz fibrations.
\end{abstract}

\maketitle

\vspace{-0.2in} 
\section{Introduction}

The seminal work of Donaldson regarding Lefschetz pencils on symplectic manifolds together with Gompf's generalization of Thurston's construction of symplectic structures on surface bundles over surfaces show that closed orientable $4$-manifolds which admit Lefschetz fibrations over the $2$-sphere are precisely the closed symplectic $4$-manifolds, up to blow-ups. In contrast, broken Lefschetz fibrations, the generalization of Lefschetz fibrations where the topology of regular fibers are allowed to change in the expense of introducing a $1$-dimensional singular set, exist on all closed smooth oriented $4$-manifolds. It is therefore natural to ask how far analogues of various results on Lefschetz fibrations extend within the class of closed smooth oriented $4$-manifolds when broken Lefschetz fibrations are considered. It is the authors' contention that constraining the topology of fibers that can appear in a broken Lefschetz fibration, and then determining which $4$-manifolds can admit such fibrations is an effective way to deal with this abundance. The classification problem we under take in our article takes this path.

Monodromy factorizations of genus one Lefschetz fibrations over the $2$-sphere, up to Hurwitz equivalences, correspond to monodromy factorizations of genus one Lefschetz fibrations over the $2$-disk with trivial global monodromy.  They 
were classified by Kas and Moishezon independently \cite{K,M}. (Also see Matsumoto's work \cite{M}.) This theorem in particular implies that the $4$-manifolds that admit non-trivial genus one Lefschetz fibrations are very restricted; namely, only elliptic surfaces $E(n)$ admit them  \cite{K,M,Mo}. We consider a generalization of this famous theorem to broken Lefschetz fibrations, aiming to classify simplified broken Lefschetz fibrations with fiber genera one and zero. These simplified broken Lefschetz fibrations are the ones where we have at most one round singular circle and connected fibers, which give the honest elliptic Lefschetz fibrations when the round singular set is empty. Fixing an isomorphism $Map(T^2) \cong SL(2, \Z)$, we consider genus one Lefschetz fibrations over the $2$-disk whose global monodromy maps to $\pm \text{id}$ times a positive power of the image of a Dehn twist along a non-separating curve in the latter group. First we generalize Kas and Moishezon's classification theorem to the monodromies of  genus one Lefschetz fibrations over the $2$-disk with these non-trivial global monodromies (Theorem~\ref{thm:factor}, Theorem~\ref{thm:factorB}, and Corollary \ref{cor:classify}). To prove this result, we use a graphical method to describe genus one Lefschetz fibrations, called a chart description, which is a suitable modification of that given in \cite{KMMW}. We then give a list of all closed smooth oriented $4$-manifolds admitting a genus one minimal simplified broken Lefschetz fibration (Theorem \ref{classifymanifolds} and Corollary \ref{classifymanifolds2}), using the handlebody descriptions studied by the first author in \cite{B1}.

\section{Preliminaries}
\subsection{Broken Lefschetz fibrations} \

Let $X$ and $\Sigma \, $ be compact connected oriented manifolds with/without boundary of dimension four and two, respectively, and $f : X \to \Sigma \, $ be a smooth surjective map with $f^{-1}(\partial \Sigma) = \partial X$. The map $f$ is said to have a Lefschetz singularity at a point $x$ contained in a discrete set $C \subset Int(X)$, if around $x$ and $f(x)$ one can choose orientation preserving charts so that $f$ conforms the complex local model 
\[(u, v) \to u^2 + v^2 \,. \]
The map $f$ is said to have a \textit{round singularity} along an embedded $1$-manifold $Z \subset Int(X) \setminus C$ if around every $z \in Z$, there are coordinates $(t, x_1, x_2, x_3)$ with $t$ a local coordinate on $Z$, in terms of which $f$ is given by 
\[(t, x_1, x_2, x_3) \to (t, x_1^2 - x_2^2 - x_3^2) \,.\]
We call the image $f(Z) \subset Int(\Sigma)$ the \textit{round image}. A \textit{broken Lefschetz fibration} is then defined as a smooth surjective map $f : X \to \Sigma \, $ which is submersion everywhere except for a finite set of points C and a finite collection of circles $Z \subset X \setminus C$, where it has Lefschetz singularities and round singularities, respectively. In particular, it is an honest surface fibration over $\partial \Sigma$. As shown in \cite{Sa, B2}, any generic map from a closed orientable $4$-manifold to the $2$-sphere can be homotoped to a broken Lefschetz fibrations over $\Sigma = S^2$, and thus, these fibrations are found in abundance. Lastly, note that whenever there is a fiber in $X$ containing a self-intersection $-1$ sphere, it can be blown-down to obtain a new broken Lefschetz fibration on $X'$ where $X= X' \# \CPb$. We will therefore focus on \textit{relatively minimal broken Lefschetz fibrations}, which do not contain such fiber components, without mentioning it any further below.

\vspace{0.2in}
\subsection{Monodromies and chart descriptions} \

Let $f : X \to D$ be a genus one Lefschetz fibration over the $2$-disk, $\Delta$ the set of critical values. Fix a base point $y_0$ in $\partial D$ so that $f^{-1}(y_0)$ is a torus whose mapping class group is used for the monodromy representation 
\[ \rho: \pi_1(D \setminus \Delta, y_0) \to Map (f^{-1}(y_0)) = Map(T^2) . \]  
Moreover we shall identify  $Map(T^2) \cong SL(2, \Z)$ as explained in Section~3. Now our monodromy representation is 
\[ \tilde \rho: \pi_1(D \setminus \Delta, y_0) \to SL(2, \Z). \]  
We denote by $\mu$ $(= \mu(f))$ the the global monodromy $\rho(\partial D)$, and by 
$\tilde \mu$  $(= \tilde \mu(f))$ the global monodromy $\tilde \rho(\partial D) $ in $SL(2, \Z)$.  

A \textit{Hurwitz arc system} for $\Delta=\{ y_1, \dots, y_n\}$ is an $n$-tuple, $(A_1, \dots, A_n)$, of embedded arcs in $D$ connecting $y_0$ and the critical values $y_1, \dots, y_n$ such that $A_i \cap A_j =\{y_0\}$ for $i \neq j$, and $A_1, \dots, A_n$ appear around $y_0$ in this order. It determines an $n$-tuple,  $(x_1, \dots, x_n)$,  of generators of $\pi_1(D \setminus \Delta, y_0)$, called a \textit{Hurwitz generator system}.  Then we call $(\rho(x_1), \dots, \rho(x_n))$ or $(\tilde\rho(x_1), \dots, \tilde\rho(x_n))$ a  \textit{Hurwitz system} of $f$, or a 
 \textit{monodromy factorization}. 
 
 We use the convention as follows: 
\[(a_1, \dots, a_n) \cdot (b_1, \dots, b_m) := (a_1, \dots, a_n, b_1, \dots, b_m) , \  \text{and} \]
\[(a_1, \dots, a_n)^m := (a_1, \dots, a_n)\cdot \cdots \cdot (a_1, \dots, a_n) \, ,\]
the concatenation of $m$ copies.

Chart description was first introduced in order to describe $2$-dimensional braids \cite{K92, K02}, and was generalized to a method describing monodromy representations of various topological objects \cite{K07}. A remarkable application of this method was a new proof of the classification of monodromies of genus one Lefschetz fibrations over the $2$-sphere \cite{KMMW}. This is equivalent to classifying monodromies of genus one Lefschetz fibrations over the $2$-disk with the trivial global monodromy. For the purpose of this paper, we need to classify  genus one Lefschetz fibrations over the $2$-disk with certain non-trivial global monodromies.  

\begin{figure}[h!] 
\begin{center}
\includegraphics[scale=0.65]{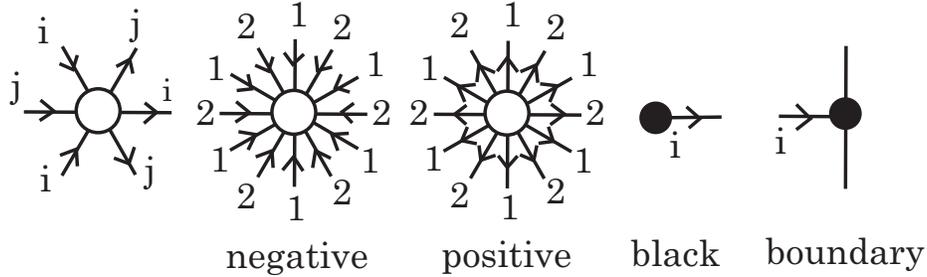}
\caption{\small Vertices of a chart; a degree-6 vertex, a negative degree-$12$ vertex, a positive degree-$12$ vertex, a black vertex, and a boundary vertex}
\label{fg_sk_1}
\end{center}
\end{figure}

\begin{defn}\label{def:chart}{\rm 
A {\it chart} is a finite graph $\Gamma$ in $D$ (possibly empty or with {\it hoops} that are closed edges with no vertices), whose edges are labeled with $1$ or $2$ and oriented so that the following conditions are satisfied:  
\begin{enumerate}
\item[(1)] The degree of each vertex is $1$, $6$ or $12$.  
\item[(2)] For a degree-six vertex $v$, the six incident edges are labeled alternately with $1$ and $2$; and three consecutive edges are oriented inward and the other three are oriented outward (see Figure~$\ref{fg_sk_1}$ where $\{ i, j \} = \{ 1, 2\}$). 
\item[(3)] For a degree-$12$ vertex $v$, the twelve incident edges are labeled alternately with $1$ and $2$; and all edges are oriented inward or all edges are oriented outward (see Figure~$\ref{fg_sk_1}$).  
\item[(4)] $\Gamma \cap \partial D$ is empty or consists of some degree-one vertices of $\Gamma$.  Moreover we assume that $\Gamma$ misses $y_0$.  
\item[(5)] For a degree-one vertex $v$ in the interior of $D$, the incident edge is oriented outward.  
\end{enumerate}
A degree-one vertex of $\Gamma$ is called a {\it black vertex} if it is in the interior of $D$, or a {\it boundary vertex} if it is on $\partial D$.  A degree-$12$ vertex is of {\it negative type} or {\it positive type} if the incident edges are oriented inward or outward, respectively. 
}\end{defn}

\begin{rk}
Definition~$\ref{def:chart}$ is slightly different from that of \cite{KMMW}.  The $4$th condition in \cite{KMMW} is that $\Gamma \cap \partial D$ is empty.  We modified it so that we can treat genus one Lefschetz fibrations with non-trivial global monodromies.  The $5$th condition is introduced here to allow only positive Dehn twists to appear in local monodromies, i.e. we only allow honest Lefschetz singularities.  When we allow negative Dehn twists, too, then the definition of a chart should be given without the $5$th condition. In this case, Proposition~$\ref{thm:chartdescription}$ and its proof are still valid.  However, Theorem~$\ref{thm:factor}$ does not work, because we can insert $(s_1, s_1^{-1})$ into any monodromy factorization.  
\end{rk}

A chart $\Gamma$ determines a homomorphism 
$\rho_\Gamma : \pi_1(D \setminus \Delta_\Gamma, y_0) \to SL(2, \Z)$, where $\Delta_\Gamma$ is the set of black vertices, as follows:  Let $\eta: [0,1] \to D \setminus \Delta_\Gamma$ be a map with $\eta(0)=\eta(1)=y_0$.  
Up to homotopy, assume that it intersects with $\Gamma$ transversely.  For each intersection, we associate a letter $s_i^\epsilon$ if the edge of $\Gamma$ at the intersection is labeled with $i$ and $\epsilon$ is $+1$ or $-1$ according to the orientation of the edge is from left to right or right to left along the direction of $\eta$.  Read these letters along $\eta$ and we obtain a word in $\{s_1, s_1^{-1}, s_2, s_2^{-1} \}$, which we call the {\it intersection word} of $\eta$ with respect to $\Gamma$ and denote it by $w_\Gamma (\eta)$.  The element of $SL(2, \Z)$ represented by this word is uniquely determined, by which we define $\rho_\Gamma ([\eta])$ (cf. \cite{K07, KMMW}).     
Here  
we regard  $s_1$ and $s_2$ as the  matrices 
$$ 
s_1 = \left( \begin{array}{cc}
1 & 0 \\
1 & 1 
\end{array} \right) 
\quad \mbox{and} \quad 
s_2 = \left( \begin{array}{cc}
1 & -1 \\
0 & 1 
\end{array} \right), 
$$
so that the group $SL(2, \Z)$ has a presentation 
$$ \langle 
s_1, s_2 \, | \, s_1s_2s_1(s_2s_1s_2)^{-1}, (s_1 s_2)^6 \rangle.$$

\begin{prop}\label{thm:chartdescription}
For any genus one Lefschetz fibration over the $2$-disk, $f: X \to D$, there exists a chart $\Gamma$ such that the monodromy representation of $f$ is equal to $\rho_\Gamma$.  
\end{prop}

\begin{proof}  This is a consequence of Theorem~$5$ of \cite{K07}.  In Theorem~$5$ and Example~$3$ of \cite{K07}, genus one Lefschetz fibrations were allowed to have singular fibers whose local monodromies were negative Dehn twists (also known as achiral Lefschetz singularities). Since we constrain local monodromies to be positive Dehn twists here, the black vertices should have incident edges oriented outward, the $5$th condition of Definition~$\ref{def:chart}$.  
\end{proof}

We call a chart $\Gamma$ as in Proposition~$\ref{thm:chartdescription}$ a 
{\it chart description} of the Lefschetz fibration $f: X \to D$.  Such a chart is not unique. There are some moves on charts, called \textit{chart moves}, that do not change the Lefschetz fibration \cite{K07}. 
In Section~\ref{section:chartmainpart}, we show that 
any chart of $f$ can be changed to a certain standard form by chart moves (Theorem~$\ref{thm:chart}$).  As an application, we will  prove the following theorem.

\begin{thm}\label{thm:factor} For a genus one Lefschetz fibration over the $2$-disk, $f: X \to D$,  with a monodromy representation $\tilde \rho: \pi_1(D \setminus \Delta, y_0) \to SL(2, \Z)$ with global monodromy  
$\tilde \mu$. Suppose that $\tilde \mu$ is $s_1^k$ or $(s_1 s_2)^3 s_1^k$ in $SL(2, \Z)$ for a non-negative integer $k$. Then after applying elementary transformations, $f$ has monodromy factorization equal to $(s_1, s_2)^{6p} \cdot (s_1)^k$ or to $(s_1,  s_2)^{6p+3} \cdot (s_1)^k$, respectively.  
\end{thm}

This theorem is equivalent to the following.   

\begin{thm}\label{thm:factorB} 
Let $k$ be a non-negative integer.  
Let $(g_1, \dots, g_n)$ be an $n$-tuple of elements of $SL(2, \Z)$ which are conjugates of $s_1$.  
\begin{itemize}
\item[(1)]  
If $g_1 \cdots g_n= s_1^k$ in $SL(2, \Z)$ for some $k$, then $p:= (n-k)/12$ is a non-negative integer and 
$(g_1, \dots, g_n)$ can be changed to  
$(s_1, s_2)^{6p} \cdot (s_1)^k$ by elementary transformations. 
\item[(2)]  
If $g_1 \cdots g_n= (s_1 s_2)^3 s_1^k$ in $SL(2, \Z)$ for some $k$, then $p:= (n-6- k)/12$ is a non-negative integer and 
$(g_1, \dots, g_n)$ can be changed to  
$(s_1,  s_2)^{6p+3} \cdot (s_1)^k$ by elementary transformations. 
\end{itemize}
\end{thm}

The case of $k=0$ in (1) of Theorem~$\ref{thm:factorB}$ is the famous theorem due to Moishezon \cite{Mo}.  
The reason why we consider that the global monodromy $\tilde \mu$ is $s_1^k$ or $(s_1 s_2)^3 s_1^k$ will be  explained in 
Section~\ref{section:SBLF} (Theorem~$\ref{SBLFmonodromy}$).

\vspace{0.2in}
\subsection{Handlebody descriptions of broken Lefschetz fibrations} \

A broken Lefschetz fibration over the $2$-disk with connected round singular set and round image an embedded curve parallel to the boundary of the $2$-disk can be depicted rather easily using handlebodies. These assumptions yield to having a Lefschetz fibration over the $2$-disk, and a round $2$-handle attached to it. Recall that a round $2$-handle is a pairwise $2$-handle attachment parametrized along $S^1$. That is, we glue an $S^1 \x D^3$ to a Lefschetz fibration over the $2$-disk in one of the two possible ways: There are two splittings of the $D^3 = D^2 \x D^1$ bundle over $S^1$ into a $D^2$-bundle and a $D^1$-bundle over $S^1$, as classified by the homotopy classes of mappings from $S^1$ into the Grassmannian $\textbf{G}(3, 2)$. Since $\pi_1(\textbf{G}(3, 2)) = \Z_2$, we get two splittings of this sort up to isotopy, each specifying a $3$-dimensional $2$-handle structure on all $D^3$ fibers of the initial (trivial) bundle $S^1 \x D^3 \to S^1$. The boundary restriction on the first component gives an $S^{1} \x D^{1}$ subbundle over $S^1$. The total space $L$ of this subbundle is a submanifold of $S^1 \x D^3$. Hence a \emph{$4$-dimensional round $2$-handle} is a copy of $S^1 \x D^{3}$, attached to the boundary of an $4$-dimensional manifold $X$ by an embedding of $L \hookrightarrow \partial X$. Round handles corresponding to the trivial splitting of the $D^3$ bundle over $S^1$ are called \emph{regular} or {untwisted} round $2$-handles, whereas those corresponding to the nontrivial splitting are called \emph{twisted}. 

Regarding the circle factor of a regular (untwisted) round $2$-handle $S^1 \x D^2 \x D^1$ as the union of a $0$-handle and a $1$-handle, we can express an untwisted round \linebreak $2$-handle as the union of a $4$-dimensional $2$-handle $H_2$ and a $3$-handle $H_3$. For a twisted round $2$-handle one obtains a similar decomposition. The splittings imply the difference: the $3$-handle goes over the $2$-handle geometrically twice and algebraically zero times in the untwisted case, and both geometrically and algebraically twice in the twisted case. (The reader can turn to \cite{B1} for the details, and for general round handles.)

Let us now describe the Kirby diagrams where one attaches a round $2$-handle to a Lefschetz fibered $4$-manifold with boundary. The round $2$-handle attachment to a surface fibration over a circle that bounds a Lefschetz fibration is realized as a fiberwise $2$-handle attachment. The attaching circle of the $2$-handle $H_2$ of a round $2$-handle is a simple closed curve $\gamma$ on a regular fiber, which is preserved under the monodromy of this fibration up to isotopy. Since this attachment comes from a fiberwise handle attachment, $H_2$ should have fiber framing zero. As usual, we do not draw the $3$-handle $H_3$ of the round $2$-handle, which is forced to be attached in a way that it completes the fiberwise $2$-handle attachments. The difference between the untwisted and twisted cases is implicit: It is distinguished by whether the curve $\gamma$ is mapped to $\gamma$ or $-\gamma$ under a self-diffeomorphism of the fiber determined by the monodromy; yielding an untwisted or a twisted round $2$-handle, respectively.

The usefulness of working with round $2$-handles is that one can depict any Lefschetz fibration over a disk together with a round $2$-handle attachment via explicit Kirby diagrams. We first draw the Lefschetz $2$-handles following the monodromy data on a regular diagram of $D^2 \x \Sigma_g$ (where $\Sigma_g$ is the regular fiber) with fiber framings $-1$, then attaches $H_2$ with fiber framing $0$ and includes an extra $3$-handle to complete it to a round $2$-handle. We draw the standard Kirby diagram where the $1$-handles depict the fiber and thus we can match the fiber framings with the blackboard framings.

\begin{figure}[ht] 
\begin{center}
\includegraphics[scale=0.67]{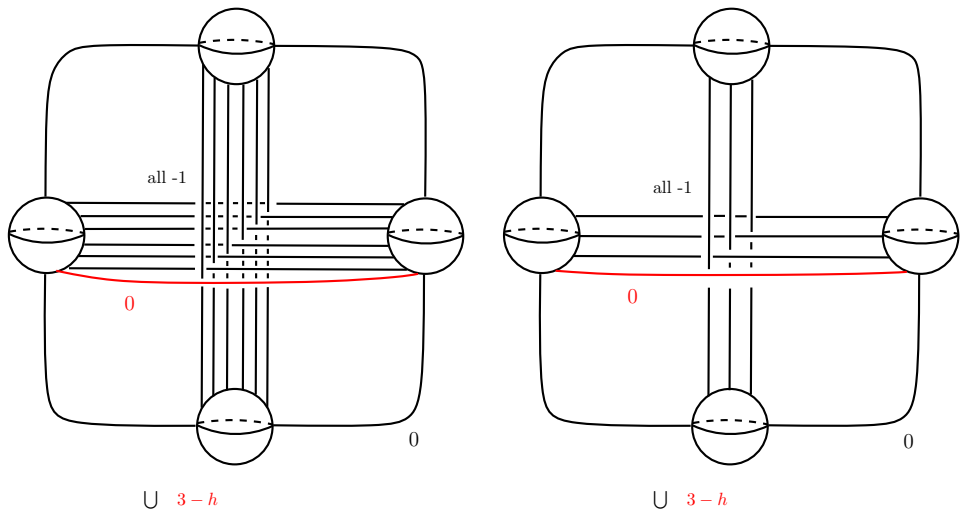}
\caption{\small Regular and twisted round $2$-handle attachments to elliptic Lefschetz fibrations over $D^2$ with monodromy factorization $(t_a, t_b)^6$ and $(t_a, t_b)^3$, respectively. Red handles make up the round $2$-handles.}
\label{round2handles}
\end{center}
\end{figure}

To illustrate our descriptions above, let us consider the two examples given in Figure \ref{round2handles}. In the first example the round $2$-handle is attached to an elliptic Lefschetz fibration with global monodromy isotopic to $\text{id}$, so $\gamma$ (given by the red $0$-framed $2$-handle) is mapped onto itself with the same orientation. Therefore it is a regular round $2$-handle. Whereas in the second example the global monodromy is isotopic to $-\text{id}$, mapping $\gamma$ to $-\gamma$. Thus, this is a twisted round $2$-handle attachment. Both of these examples will be revisited later in the paper.

\vspace{0.2in}
\section{Simplified broken Lefschetz fibrations } \label{section:SBLF}

Let $X$ be a closed orientable $4$-manifold. We will put further constraints on the broken Lefschetz fibrations in consideration to have a more tractable family. First, we ask the round singular set to be connected, i.e. to consist of one circle only, and its image on $S^2$ to be embedded. Second, we ask all the regular fibers to be connected. (Note that in general round singularities can give rise to disconnected regular fibers.) This gives a decomposition of the broken Lefschetz fibration into three pieces; a genus $g-1$ Lefschetz fibration over a $2$-disk we call the \emph{lower side}, a genus $g$ Lefschetz fibration over a $2$-disk called the \emph{higher side}, and a round cobordism between them containing the round singular set in the interior, where $g > 0$. The third, and the last condition we impose is to have all the Lefschetz singularities on the higher side, which equivalently means that the lower side consists of a trivial genus $g-1$ surface bundle. Broken Lefschetz fibrations satisfying these properties were extensively studied in \cite{B1}, under the name \emph{simplified broken Lefschetz fibrations} (abbreviated SBLF), which we will adapt herein as well. (Careful reader will notice that in \cite{B1}, the connectivity of fibers was not built into the definition of a simplified fibration. It was later shown in \cite{B2} that this could always be achieved after a homotopy.) Moreover, if the highest genus of a regular fiber in a given SBLF is $g$, we will call it a \emph{genus $g$ simplified broken Lefschetz fibration.} Observe that a genus $g$ SBLF can possibly have empty round singular set, in which case it is an honest genus $g$ Lefschetz fibration over $S^2$. 

Relying on the work of Gay-Kirby in \cite{GK}, one can always obtain \emph{achiral} broken Lefschetz fibrations over $S^2$ with embedded round image. These achiral Lefschetz singularities can then be replaced locally by broken Lefschetz fibrations, as argued in \cite{B3, L}. (Or alternatively the construction in \cite{AK} can be employed, where achiral singularities are already avoided.) Furthermore, as observed by Williams, one can homotope such a broken Lefschetz fibration to one with connected round singular set and embedded round image \cite{W}. Next, the \textit{flip-and-slip} move of \cite{B2} can be employed to obtain a homotopic broken Lefschetz fibration with the same properties and with only connected fibers. Finally, one can \textit{push} the Lefschetz singularities to the higher side, as argued in \cite{B1}. In short, \textit{there is always a simplified broken Lefschetz fibration on any $X$}.

The monodromy representations of SBLFs are simple. Let $Map_{\gamma}(\Sigma_g)$ be the subgroup of $Map(\Sigma_g)$ that consists of elements which fix the embedded curve $\gamma$, up to isotopy. Then there is a natural homomorphism 
\[ \phi_{\gamma}: Map_{\gamma}(\Sigma_g) \to Map(\Sigma_{g-1}) \, .\]
Observe that our assumption on the connectivity of fibers implies that $\gamma$ is a nonseparating curve. Define $S_g$ to be the set of pairs $(\mu, \gamma)$ such that $\mu \in Map_{\gamma}(\Sigma_g)$ and $\mu \in Ker \, (\phi_{\gamma})$. Recall that when the fiber genus is at least two, fiber-preserving gluing maps are determined uniquely up to isotopy. Hence, given any tuple \linebreak $(\mu, \gamma) \in S= \bigcup_{g \geq 3} S_g$, we can construct a unique SBLF. If $g < 3$, then one also needs to tell how the pieces are glued along the low genera surface bundles over circles; amounting to two possible choices for $g=0$ pieces and $\Z^2$-choices for $g=1$ pieces. 

The map $\phi_{\gamma}: Map_{\gamma}(\Sigma_g) \to Map(\Sigma_{g-1})$ above factors as 
\[ \psi_{\gamma}:  Map_{\gamma}(\Sigma_g) \to Map(\Sigma_g \setminus N) \ \ \text{and} \ \ \ \varphi_{\gamma}: Map(\Sigma_g \setminus N) \to Map(\Sigma_{g-1}), \] 
where $N$ is an open tubular neighborhood of $\gamma$ away from the other vanishing cycles. (The middle group does not need to fix the boundaries.) It is easy to see that the map $\psi_{\gamma}$ has kernel isomorphic to $\Z$. When we have a SBLF, the kernel of $\varphi_{\gamma}$ is isomorphic to the braid group on $\Sigma_{g-1}$ with $2$-strands, by definition. 

Now let us assume that $f: X \to S^2$ is a SBLF with $F_g = T^2$. We fix two generators $a$ and $b$ of $\pi_1(T^2) \cong \Z^2$ and an isomorphism $Map(T^2) \cong SL(2, \Z)$, such that the positive Dehn twist $t_a$ is mapped to $s_1$ and $t_b$ to $s_2$, where
$$ 
s_1 = \left( \begin{array}{cc}
1 & 0 \\
1 & 1 
\end{array} \right) 
\quad \mbox{and} \quad 
s_2 = \left( \begin{array}{cc}
1 & -1 \\
0 & \, 1 
\end{array} \right).  
$$
What underlies our choices here is the convention of \cite{KMMW}. So the curves $a$ and $b$ correspond to $(0,1) ^{\, T}$ and $(1,0) ^{\, T}$. Let $\tilde{\mu}$ be the image of the globaly monodromy $\mu$ of the Lefschetz fibration on the higher side of $f: X \to S^2$ under the chosen isomorphism $Map(T^2) \cong SL(2, \Z)$. 

Without loss of generality, we can assume that the non-seperating curve $\gamma$ is equal to $a$, so the above condition translates to having 
\[\tilde{\mu} \, (0,1)^{\, T} = (0,1)^{\, T} \ \text{or} \ \  \tilde{\mu} \, (0,1)^{\, T} = -(0,1)^{\, T} .\]
For these two cases, we respectively get:
$$ 
\tilde{\mu} = \left( \begin{array}{cc}
1 & 0\\
m &  1 
\end{array} \right) 
\quad \mbox{or} \quad 
\left( \begin{array}{cc}
-1 & \, 0 \\
\, n & -1 
\end{array} \right),  
$$
where $m,n$ are arbitrary integers. Note that the former corresponds to having a regular round handle cobordism, whereas the latter amounts to a twisted one. (See \cite{B1}.) 

\begin{figure}[ht]
\begin{center}
\includegraphics[scale=1.2]{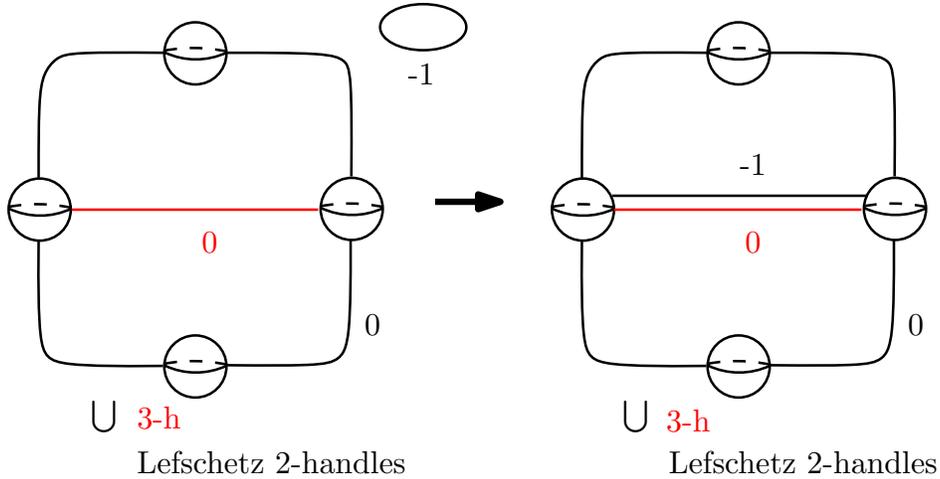}
\caption{\small Blow-up of a broken Lefschetz fibration over a $2$-disk.}
\label{blowup}
\end{center}
\end{figure}

If the diagonal entries of $\tilde{\mu}$ are $+1$, then $\tilde{\mu} = s_1^m$. If they are $-1$, then we can express $\tilde{\mu}$ as $-\text{id} \, s_1^{-n}= (s_1 s_2)^3 s_1^{-n}$. For $m \geq 0$ (resp. $n \leq 0$), the right hand side of the first (resp. second) expression corresponds to a product of right-handed Dehn twists, and in turn, to a monodromy of genus one Lefschetz fibration over the $2$-disk. On the other hand, if $m < 0$, we can employ the following trick: Include $|m|$ right-handed Dehn twists to the factorization, which provides us with a SBLF $f': X' \# |m| \CPb \to S^2$, again with higher genus one. (See the \textit{isotropic blow-up} example of \cite{ADK}, and \cite{B1} for the handlebody argument we reproduce here.) This is due to the fact that, the introduction of each extra right-handed Dehn twist along $a$ can be seen as in Figure \ref{blowup}. Since the $0$-framed $2$-handle (drawn in red in the figure) of the $2$-handle is attached fiberwise, it does not link with any one of the Lefschetz $2$-handles, explaining why we can perform this modification without interfering with the attachment of these $2$-handles. Sliding the blow-up curve over the $0$-framed $2$-handle of the round $2$-handle, we realize it as a new Lefschetz $2$-handle, attached along the same curve $a$. The induced monodromy on the lower side does not change, and therefore it can be glued to this broken Lefschetz fibration over the $2$-disk in the same way it was done for $f: X \to S^2$. Hence, we can replace $\tilde{\mu}$ with $\tilde{\mu}' = \text{id}$, after passing to a blow-up of $X$. The very same line of arguments work for $n > 0$ case as well, where we end up replacing $\tilde{\mu}$ with $\tilde{\mu}'= -\text{id}$. (Note that $s_1^{\pm 1}$ commutes both with $-\text{id}=(s_1 s_2)^3$ and $s_1^m$, so there is no order issue.) 

We have proved:


\begin{thm} \label{SBLFmonodromy}
Let $X$ admit a genus one simplified broken Lefschetz fibration. Then, possibly after blowing-up $X$, we get a genus one simplified broken Lefschetz fibration whose global monodromy on the higher side maps to $\tilde{\mu} = s_1^k$ or to $(s_1 s_2)^3 s_1^k$ in $SL(2, \Z)$, for $k$ a non-negative integer. 
\end{thm}

\vspace{0.2in}
\section{Monodromies of genus one Lefschetz fibtations over a $2$-disk}\label{section:chartmainpart}

Let $f : X \to D$ be a genus one Lefschetz fibration over a $2$-disk and 
\[\rho: \pi_1(D \setminus \Delta, y_0) \to Map(T^2)\] 
be its monodromy representation, where $\Delta$ is the set of critical values of $f$, $y_0 \in \partial D$ is a fixed base point. For a given isomorphism $Map(T^2) \cong SL(2, \Z)$, we get a monodromy representation 
$\tilde \rho: \pi_1(D \setminus \Delta, y_0) \to SL(2, \Z)$.  
 
By Theorem~$\ref{SBLFmonodromy}$, we may assume that the global monodromy 
$\tilde \mu$ is equal to  $s_1^k$ or to $(s_1 s_2)^3 s_1^k$ for a non-negative integer $k$.   
For simplicity, we denote by $q$ (or $q(f)$) the integer $0$ or $1$ such that $\tilde \mu = (s_1 s_2)^{3 q} s_1^k$.  

The number of critical values of $f$ is denoted by $c(f)$.    

\begin{thm}\label{thm:chart}
Let $f : X \to D$ be a genus one Lefschetz fibration with $\tilde \mu = (s_1 s_2)^{3 q} s_1^{k}$, where 
$q \in \{0,1\}$ and $k$ is a non-negative integer.  
Let $\Gamma$ be a chart description of $f$.   
\begin{itemize}
\item[(1)] 
By chart moves, the chart $\Gamma$ can be changed to a chart written as 
$N^{p} \amalg (U_1U_2)^{3 q} U_1^{k}$ for some non-negative integer $p$.   
\item[(2)] 
In $(1)$, the number $p$ is uniquely determined and it is equal to \\ $(  c(f) - 6 q - k)/12$. 
\end{itemize}
\end{thm}

\begin{figure}[ht] 
\begin{center}
\includegraphics[scale=0.6]{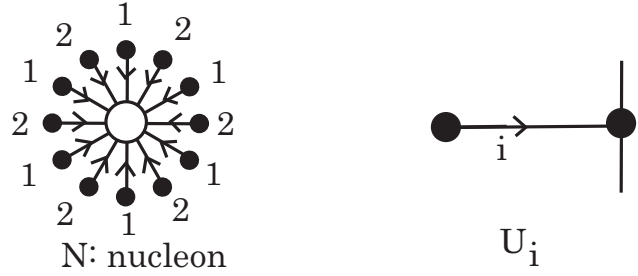}
\caption{\small Left: a nucleon,  Right: $U_i$ $(i=1,2)$}
\label{fg_sk_2}
\end{center}
\end{figure}

\begin{figure}[ht] 
\begin{center}
\includegraphics[scale=0.5]{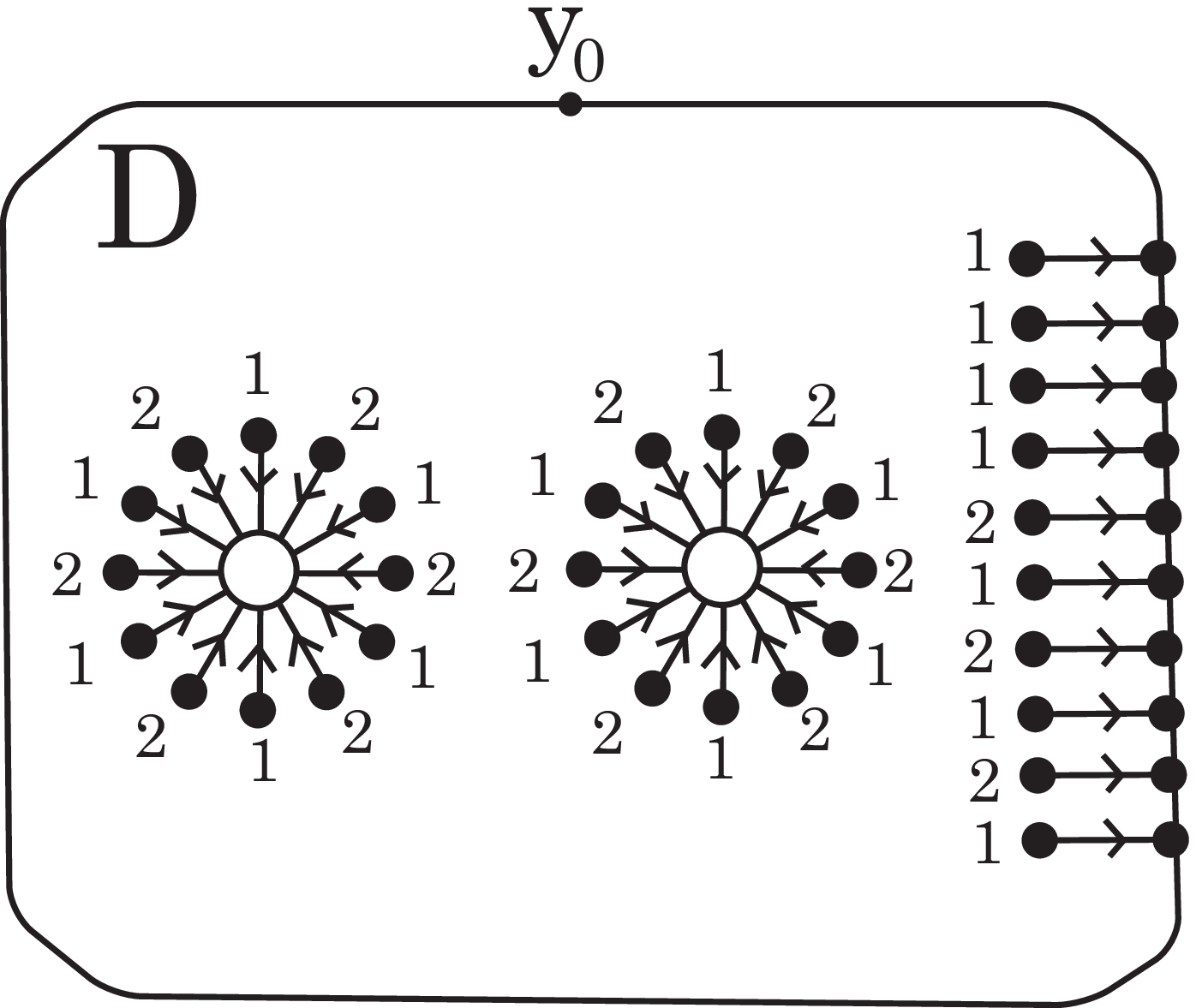}
\caption{\small $N^{2} \amalg (U_1U_2)^{3} U_1^{4}$}
\label{fg_sk_3}
\end{center}
\end{figure}

Here $N$ is a chart consisting of a single degree-$12$ vertex of negative type and 12 black vertices together with 12 edges, which we call a {\it nucleon}, and $U_i$ $(i=1,2)$  is a chart consisting of a black vertex and a boundary vertex with an edge labeled with $i$ (see Figure~$\ref{fg_sk_2}$).  The chart $(U_1U_2)^{3 q} U_1^{k}$ is the union of  some copies of $U_1$ and $U_2$ appearing along $\partial D$ in this order.  For example, $N^{2} \amalg (U_1U_2)^{3} U_1^{4}$ is as in Figure~$\ref{fg_sk_3}$. 

As a corollary to Theorem~$\ref{thm:chart}$, we obtain Theorem~$\ref{thm:factor}$:  

\noindent {\it Proof of Theorem~$\ref{thm:factor}$.}  Let $f$ be a genus one Lefschetz fibration with $\tilde \mu = (s_1 s_2)^{3 q} s_1^{k}$, where 
$q \in \{0,1\}$ and $k$ is a non-negative integer.  Take a chart description $\Gamma$ of $f$.  By Theorem~$\ref{thm:chart}$, we may assume that $\Gamma$ is  $N^{p} \amalg (U_1U_2)^{3 q} U_1^{k}$.  Taking a Hurwitz generator system of 
$\pi_1(D \setminus \Delta_\Gamma, y_0)$ in an obvious way, $N^p$ yields a factorization $(s_1, s_2)^{6p}$ and 
$(U_1U_2)^{3 q} U_1^{k}$ yields $(s_1, s_2)^{3 q} \cdot (s_1)^k$.  \qed 

\begin{cor}\label{cor:classify}
Let $f$ and $f'$ be genus one Lefschetz fibration over a $2$-disk with 
\[\tilde \mu (f) = \tilde \mu (f') = (s_1 s_2)^{3 q} s_1^{k}\] 
for $q \in \{0,1\}$ and a non-negative integer $k$. Then $f$ and $f'$ are equivalent if and only if $c(f)= c(f')$.  
\end{cor}

The remainder of this section is devoted to proving Theorem~$\ref{thm:chart}$.   

\begin{lem} \label{lem:subword}
Let $W$ be a word $(s_1 s_2)^{3q} s_1^k$ with $q \in \{0,1\}$ and $k$ a non-negative integer.  
Let $W'$ be a subword of $W$.  If $W'= 1$ in $SL(2, \Z)$, then $W'$ is the empty word.  
\end{lem}

\begin{proof} By a direct calculation, we see that if a subword $W''$ of $(s_1 s_2)^3$ is equal to $s_1^n$ for some $n\in \Z$ then $W''$ is $s_1^m$, for $m \in \{0,1,2,3\}$.  Thus $W'$ should be the empty word.
\end{proof} 

An edge of a chart is said to be of type {\it $(1,6)$}, {\it $(1,12)$} or {$(1,\partial)$} if the source is a black vertex and the target is a degree-six vertex, a degree-$12$ vertex or a boundary vertex, respectively.  An edge of a chart is said to be of type {\it $(n,6)$}, {\it $(n,12)$} or {$(n,\partial)$} where $n \in \{6, 12\}$ if the source is a degree-$n$ vertex and the target is a degree-six vertex, a degree-$12$ vertex or a boundary vertex, respectively.  

We can now prove:
 
\vspace{0.2cm} 
\noindent {\it Proof of Theorem~$\ref{thm:chart}$.}   
First, by chart moves of type $\partial$ defined in \cite{K07} (Fig. 13), we can change $\Gamma$ so that $\partial \Gamma = \partial ((U_1U_2)^{3 q} U_1^{k})$. This is possible because $\rho_\Gamma (\partial D) = \tilde \mu = (s_1 s_2)^{3 q} s_1^{k}$ (see \cite{K07}). Then let $p$ be the number of degree-$12$ vertices of negative type minus the number of those of positive type. Comparing the numbers of sources and targets of the edges, we see that $c(f)=12 p + 6 q +k$. Therefore $p=(c(f)-6q-k)/\,12$ for any chart description $\Gamma'$ with $\partial \Gamma' = \partial ((U_1U_2)^{3 q} U_1^{k})$. This implies the second assertion of the theorem. 

By the argument in the proof of Theorem~21 of \cite{KMMW}, we can change $\Gamma$ by chart moves such that  the degree-$12$ vertices are all positive or all negative. 
Using the assumption that $k$ is non-negative, we can apply an argument similar to that of Lemma~22 of \cite{KMMW} to remove all edges of type $(1,6)$.   
Now every black vertex is of type $(1,12)$ or $(1, \partial)$.  

Suppose that $p$ is non-negative, i.e., there are no degree-$12$ vertices or there are $p$ degree-$12$ vertices of negative types. In this case, we assert that there are no edges of type $(6, \partial)$ or $(12, \partial)$:  Since there are no degree-$12$ vertices of positive type, there are no edges of type $(12, \partial)$. Suppose that there is an edge $e$ of type $(6, \partial)$ whose target is a boundary vertex $v_0$ and the source is a degree-six vertex $v_1$.  The three edges incident to $v_1$ oriented toward $v_1$ are edges of type $(6,6)$. Let $v_2, v_3, v_4$ be the degree-six vertices of the sources. There might be duplication in $v_1, v_2, v_3, v_4$, but at least one of $v_2, v_3, v_4$ is not $v_1$.  Continue this argument and we obtain a strictly increasing family of degree-six vertices.  Since $\Gamma$ has a finite number of degree-six vertices, this yields a contradiction. Thus there are no edges of type $(6, \partial)$ or $(12, \partial)$.  
Now all boundary vertices are targets of edges of type $(1, \partial)$, and they forms the latter part $(U_1U_2)^{3 q} U_1^{k}$. By the proof of Theorem~21 of \cite{KMMW} again, we can change the remainder into a union of nucleons. Now we have $N^{p} \amalg (U_1U_2)^{3 q} U_1^{k}$.  

Suppose that $p$ is a negative integer, i.e., all degree-$12$ vertices are of positive type. We will show that this case never happens. All black vertices are sources of edges of type $(1, \partial)$. Let $M$ be the boundary vertices of $\Gamma$ which are the targets of the edges of type $(1, \partial)$. Let $W'$ a subword of the word $(s_1 s_2)^{3 q} s_1^{k}$ obtained by deleting letters that correspond to the points of $M$. Then $W' = 1$ in $SL(2, \Z)$. (This is seen as follows. Let $C$ be a simple loop in $D$ which is close and parallel to $\partial D$.   Shrinking the $(1, \partial)$-edges, we may assume that they are between the loop $C$ and $\partial D$. Since there is no black vertices inside of $C$, $\rho_\Gamma(C) =1$ in $SL(2, \Z)$.) By Lemma~$\ref{lem:subword}$, $W'$ must be the empty word.  
Thus $\Gamma$ is chart move equivalent to $\Gamma' \amalg (U_1 U_2)^{3q} U_1^k$ for some chart $\Gamma'$ without black vertices such that $\Gamma' \cap \partial D = \emptyset$.  By a chart move, we can remove $\Gamma'$ to obtain $(U_1 U_2)^{3q} U_1^k$.  Then $p=0$, a contradiction. 

Hence we see that $\Gamma$ is chart move equivalent to $N^{p} \amalg (U_1U_2)^{3 q} U_1^{k}$.  \qed

\newpage  
\section{Four-manifolds admitting genus one SBLFs} 

We are ready to identify the total spaces of higher genus one relatively minimal SBLFs. We are going to start with the simplest cases when either the round singular set or the set of Lefschetz critical points is empty:

\noindent \textbf{No round singularity.} If in addition there are no Lefschetz critical points, then the fibration can be trivialized over each hemisphere of the base $S^2$. Thus, the total space and the fibration are obtained by gluing two copies of $T^2 \x D^2$ equipped with projection maps onto $D^2$ via some fiber-preserving diffeomorphism on their boundaries. This is equivalent to performing a multiplicity $\pm 1$ logarithmic transform along a regular fiber $T$ of the standard fibration on $T^2 \x S^2$. 

\begin{figure}[h]
\begin{center}
\hspace{-1cm}
\includegraphics[scale=0.7]{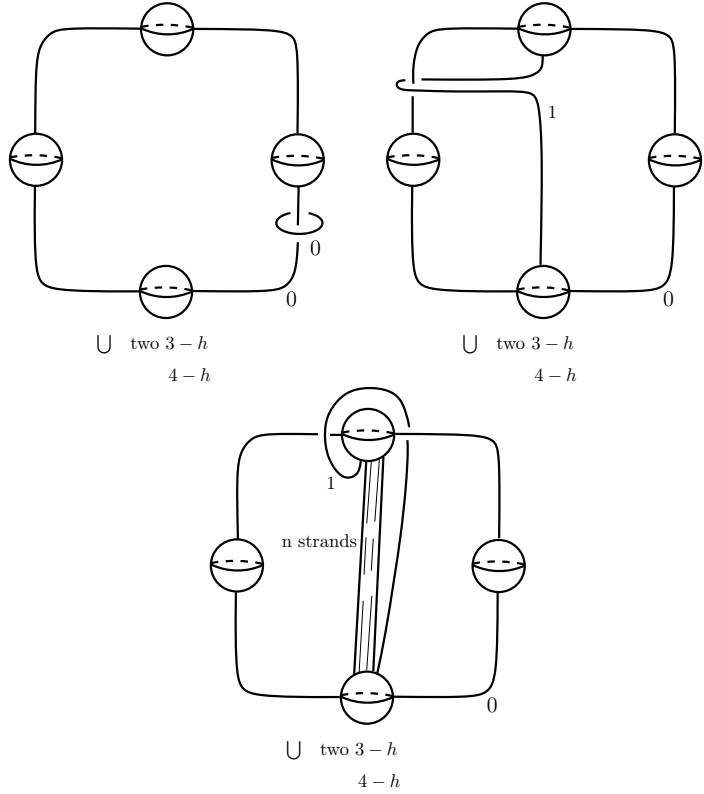}
\caption{\small Top left: The trivial fibration on $T^2 \x S^2$. Top right: The fibration on $S^1 \x S^3$ obtained from the Hopf fibration. Bottom: The locally trivial torus fibration on $S^1 \x L(n,1)$, for $n>1$.} 
\label{genusoneLFs}
\end{center}
\end{figure}

Recall that the fiber framing prescribes an isomorphism 
\[H_1(\partial \nu T)  \cong  H_1(T ; \Z) \oplus \Z \, , \] 
where the $\Z$ component is generated by the positively oriented meridian $m_T$ of $T$, so the image of $[\partial D^2]$ under the homomorphism induced by the boundary diffeomorphism is of the form $r [C] \pm [m_T]$ for some primitive curve $C$ on $T$. Note that for $r=0$ we get back the trivial fibration on $T^2 \x S^2$, and for $r=1$, we get the standard fibration on $S^1 \x S^3$ derived from the Hopf fibration on $S^3$. We can map $C$ to any primitive curve using a self-diffeomorphism of the fiber which clearly extends over $T^2 \x D^2$. A handlebody description of the total space is therefore obtained by adding a $2$-handle to the standard diagram for $T^2 \x D^2$ as shown in the last diagram given in Figure \ref{genusoneLFs}. It is easy to see that the total space in this case is $S^1 \x L(n,1)$ for some $n>1$ \cite{GS}. Finally, our choice of orientation for $L(n,1)$ was unimportant, since $S^1 \x L(n,1)$ admits an orientation preserving self-diffeomorphism composed of an orientation reversing diffeomorphism on both components.

\begin{lem}
The only closed oriented $4$-manifolds admitting a locally trivial torus fibration over $S^2$ are $S^2 \x T^2$, $S^1 \x S^3$, and $S^1 \x L(n,1)$.
\end{lem}
 
\begin{figure}[ht]
\begin{center}
\hspace{-1cm}
\includegraphics[scale=0.7]{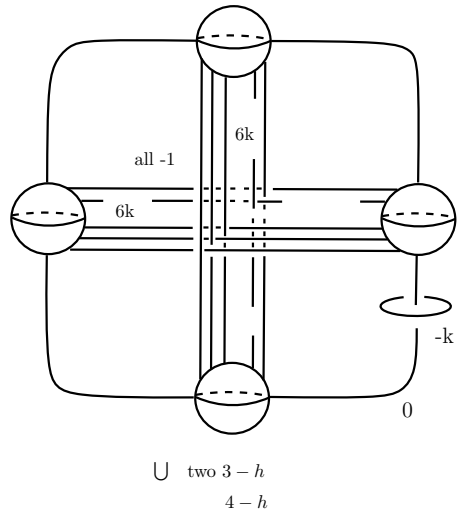}
\caption{\small Elliptic Lefschetz fibration on $E(k)$.} 
\label{genusoneLFs2}
\end{center}
\end{figure}

When there are Lefschetz singularities, we get a genus one Lefschetz fibration with monodromy factorization $(s_1, s_2)^{6k}$ in $SL(2, \Z)$, whose total spaces is the elliptic surface $E(k)$, given in Figure \ref{genusoneLFs2}. (This is the classical result of Kas, Moishezon, Matsumoto \cite{K, Mo, M}.) 

\begin{figure}[h]
\begin{center}
\includegraphics[scale=0.7]{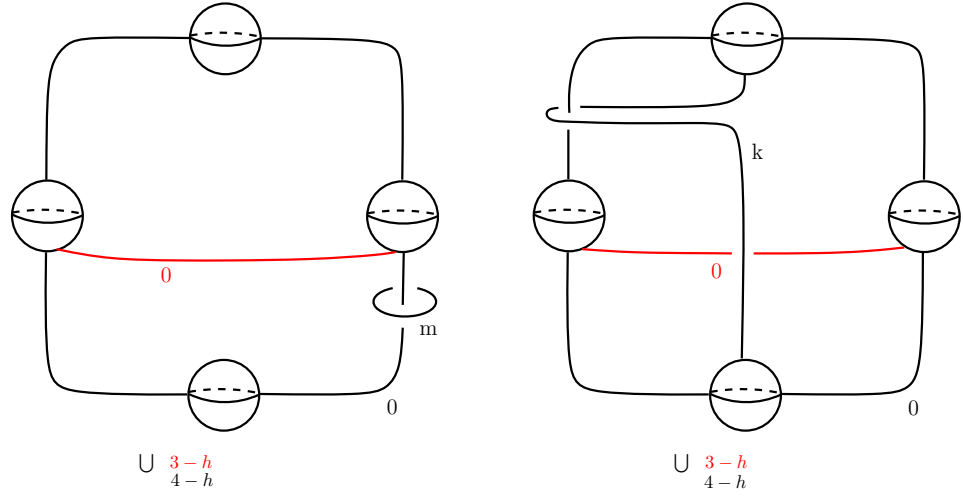}
\caption{\small Left: Total space is $S^2 \x S^2 \# S^1 \x S^3$ or $\CP \# \CPb \# S^1 \x S^3$, depending on whether $m$ is even or odd, respectively. Right: SBLF on the $4$-sphere.}
\label{trivialgenusoneSBLFs}
\end{center}
\end{figure}

\noindent \textbf{One round singular circle, no Lefschetz singularity.} Three possibilities for the total spaces in this case are $S^2 \x S^2 \# S^1 \x S^3$, $\CP \# \CPb \# S^1 \x S^3$, or $S^4$, as explored in \cite{ADK} and depicted by the Kirby diagrams in Figure \ref{trivialgenusoneSBLFs}. We will refer to these as ``standard'' broken fibrations on the corresponding $4$-manifolds. The calculus to verify the total spaces can be found in \cite{B1}.

\begin{figure}[h!]
\begin{center}
\includegraphics[scale=0.7]{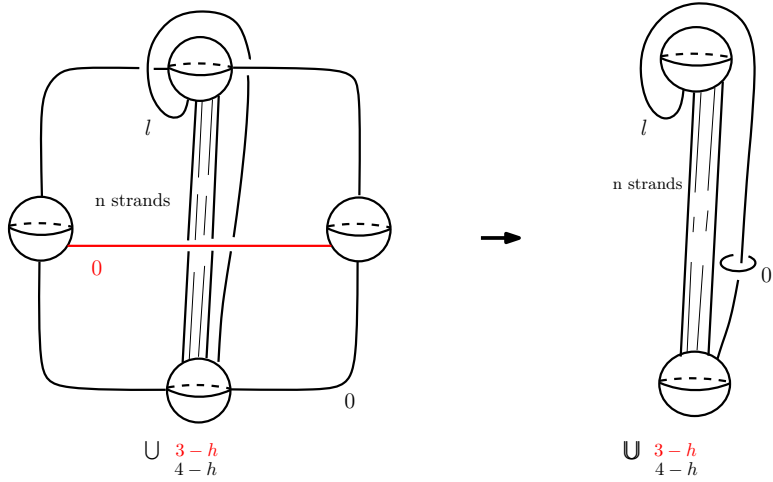}
\caption{\small Other genus one SBLFs with no Lefschetz singularities. On the right: The manifolds $L_n$ and $L'_n$ for $l$ even and odd, respectively.}
\label{trivialgenusoneSBLFs2}
\end{center}
\end{figure}

We will show that there are two more possibilities for the total spaces of such fibrations, $L_n$ and $L'_n$, completing the list in this case. To see this, observe that the handle diagram for the total space consists of the standard diagram for $T^2 \x D^2$ together with the round $2$-handle, and the $2$-handle from the lower side pulled back to this diagram via a fiber preserving self-diffeomorphism of the $T^3$ boundary. This $2$-handle can be unlinked from the $1$-handle that the $0$-framed $2$-handle of the round $2$-handle is linking with, giving us the handle diagram on the left hand side of the Figure \ref{trivialgenusoneSBLFs2}. Note that the framing of this $2$-handle depends on how the $2$-handle from the lower side is pulled back to this diagram and can attain arbitrary values $l$. We can then simplify this diagram by first sliding off both the $2$-handle coming from the lower side and the $2$-handle corresponding to the fiber over the $0$-framed $2$-handle of the round $2$-handle, and then canceling the $1$-handle that links with this $0$-framed $2$-handle. This results in the diagram on the right hand side of the Figure \ref{trivialgenusoneSBLFs2}. Note that the framing $l$ of the remaining $2$-handle linking with the $1$-handle can be made $0$ or $1$ by using the $0$-framed $2$-handle linking it once. Hence for each $n$, the number of times the $2$-handle runs over the $1$-handle, there are only two types of $4$-manifolds depending on the parity of $l$. It turns out that these $4$-manifolds $L_n$ and $L'_n$ are in fact the ones introduced by Pao \cite{P}, for $l$ even and odd, respectively, whose Kirby diagrams  described in \cite{H} match with ours.

Thus we have shown:

\begin{lem}
The only closed oriented $4$-manifolds admitting a genus one broken Lefschetz fibration with no Lefschetz singularities are $S^4$, $S^2 \x S^2 \# S^1 \x S^3$, \linebreak $\CP \# \CPb \# S^1 \x S^3$, $L_n$, and $L'_n$. 
\end{lem}

\begin{rk}
One can alternatively produce all these SBLFs from the standard one on $S^4$ by performing various log $1$ transforms along a torus fiber on the higher side and Gluck twists along an $S^2$ fiber on the lower side while preserving the fibration structure. The former yields $S^2 \x S^2 \# S^1 \x S^3$ or $L_n$, whereas the latter alters the parity and hands us $\CP \# \CPb \# S^1 \x S^3$ or $L'_n$, respectively. This can be easily verified by drawing the handle diagram of the lower side following \cite{B3} for the two cases, and then analyzing all possible fiber preserving gluings of $T^2 \x D^2$.   
\end{rk}

\noindent \textbf{Non-trivial cases.} Now we assume that neither the round singular set nor the set of Lefschetz critical points is empty. From our results in the previous sections (namely, Theorems \ref{thm:factor} and \ref{SBLFmonodromy}) it follows that, possibly after passing to a blow-up, the factorization of the global monodromy on the higher side is Hurwitz equivalent to $(t_a,  t_b)^{3n} \cdot (t_a)^k$, for some non-negative integers $n$ and $k$. 
Using the blow-up argument we gave in Section \ref{section:SBLF}, it is easy to see that when $n=0$, we in fact get the blow-ups of any the genus one SBLFs with one round singular circle covered in the previous case. Below, we will focus on the case $(t_a,  t_b)^{3n+3} \cdot (t_a)^k$ for $n$ non-negative.

\begin{figure}[ht]
\begin{center}
\includegraphics[scale=0.8]{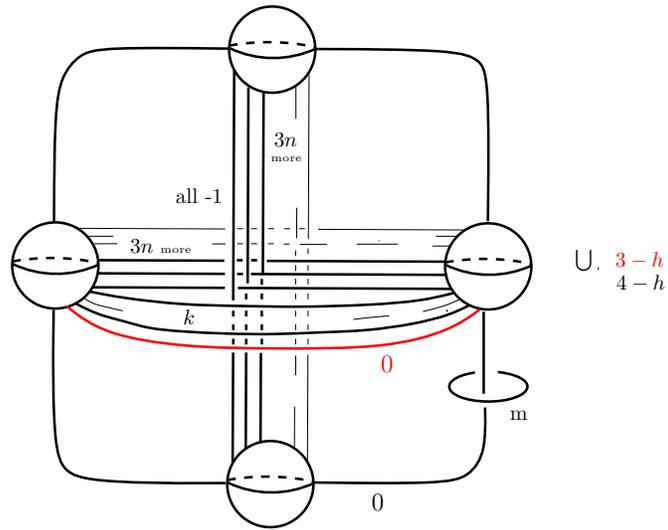}
\caption{\small The SBLF with higher side monodromy $(t_a, t_b)^{3n+3} \cdot (t_a)^k$ and a section.}
\label{calculus1a}
\end{center}
\end{figure}

Recall that to identify the total space of such a fibration, we also need to know the identification of the boundaries of the higher and lower side with the ends of the round cobordism in between. Let us first assume that the higher side is identified using the identity, and postpone the discussion of the `twisted' cases. Then there is a section of the genus one fibration on the higher side which extends through the round cobordism. This can be matched with a section of the trivial genus zero fibration on the lower side to get a global section $S$ of the broken Lefschetz fibration, with self-intersection $m$ equals the sum of the self-intersection of the disk section of the genus one fibration and the disk section of the genus zero fibration contained in $S$. The latter self-intersection can be any integer; this can be thought of the number of times the fibers are fully rotated when identifying the $S^2$ bundles over $S^1$ on the boundaries of the lower side and the round cobordism. Nevertheless, we are going to see that the total space is independent of $m$. Note that a slight extension of the classical observation of Moishezon shows that any orientation and fiber preserving self-diffeomorphism of the boundary of the higher side with global monodromy $(t_a,  t_b)^{3n}$ (and therefore of $(t_a,  t_b)^{3n} \cdot (t_a)^k$) can be extended to the interior (an elementary proof of which is given in \cite{K}), so we can keep assuming that the $2$-handle of the round $2$-handle is attached as shown in Figure \ref{calculus1a}.

Figure \ref{calculus1a} gives a handlebody description of such a SBLF. The section $S$ is represented by the $2$-handle with framing $m$, linking the $0$-framed $2$-handle capping off the obvious genus one surface given by the two $1$-handles. Observe that there are two types of Lefschetz handles; $3m$ pairs of Lefschetz handles with vanishing cycles along $a$ and $b$ attached in an alternating fashion with the given linking pattern (for $(t_a, t_b)^{3n+3}$), and $k$ Lefschetz handles attached along the curve $a$ afterward (for $(t_a)^k$). The first three pairs are drawn thicker in the figure. We will demonstrate our Kirby calculus arguments first using these three pairs, before handling the remaining $3n$ pairs inductively. Lastly, the round $2$-handle is composed of the $0$-framed $2$-handle attached along $a$ and the $3$-handle, both of which are given in red in the figure. 

\begin{figure}[h!]
\begin{center}
\includegraphics[scale=0.8]{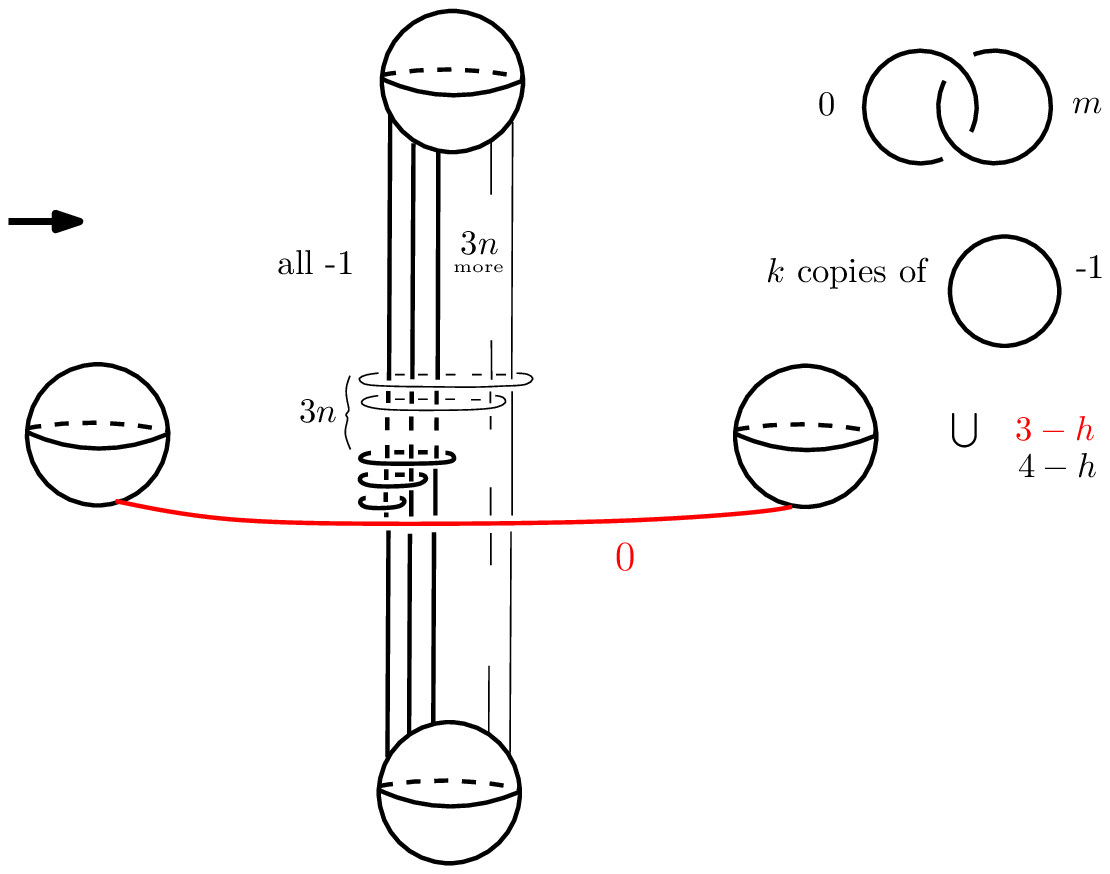}
\caption{\small }
\label{calculus1b}
\end{center}
\end{figure}

In order to arrive to a simpler handlebody diagram of the closed orientable \linebreak $4$-manifold admitting this SBLF, we start with sliding-off all the $2$-handles that were linking the $1$-handle carrying the curve $a$, using the $0$-framed $2$-handle of the round $2$-handle. This includes sliding the $0$-framed $2$-handle corresponding to the fiber over the $0$-framed $2$-handle of the round $2$-handle twice, and then isotoping it away from the other $1$-handle together as well (while the $m$-framed $2$-handle is dragged away with it). The resulting diagram is given in Figure \ref{calculus1b}.

\begin{figure}[htp]
\begin{center}
\includegraphics[scale=0.9]{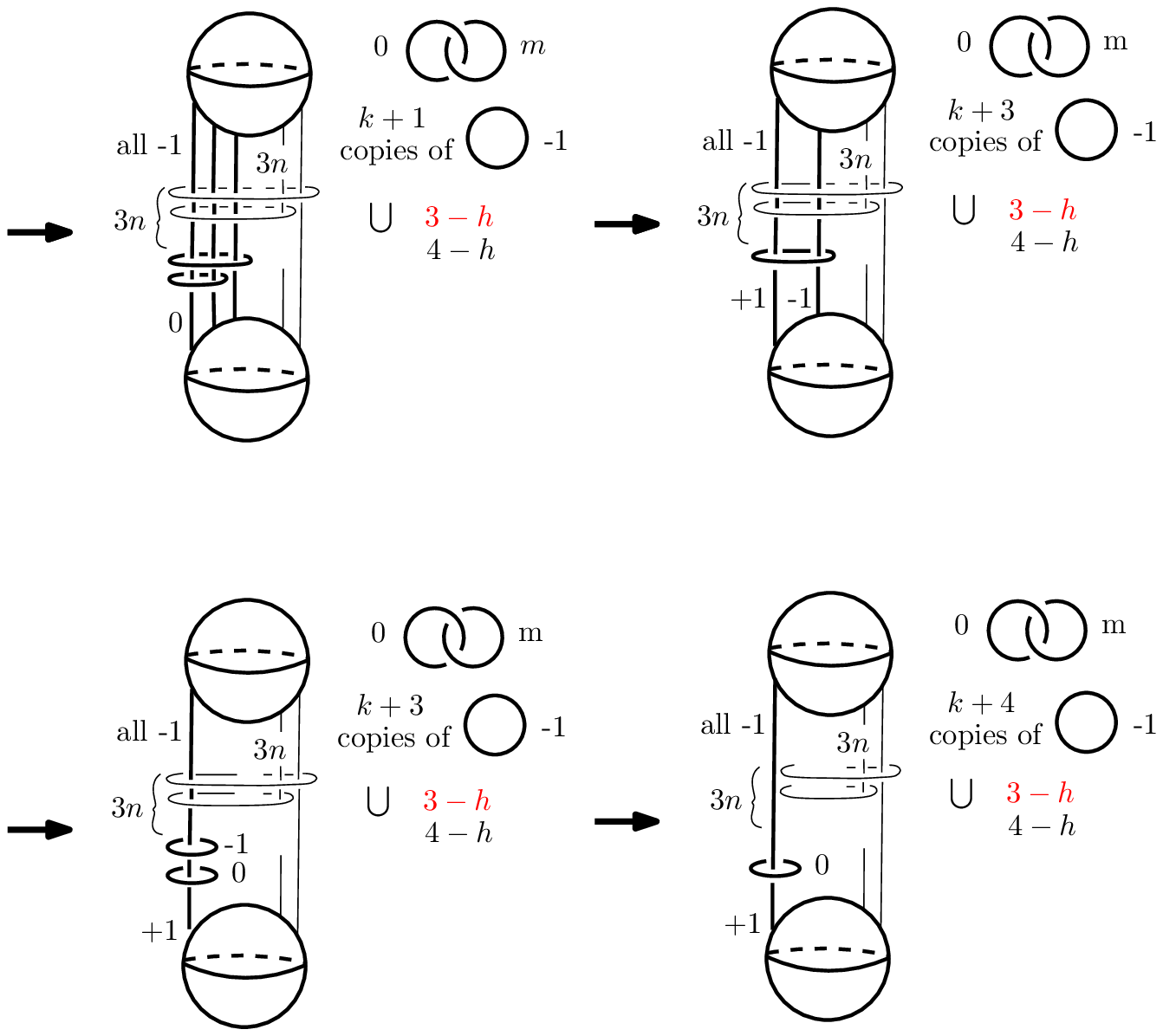}
\caption{\small }
\label{calculus2a}
\end{center}
\end{figure}

\begin{figure}[htp]
\begin{center}
\includegraphics[scale=0.9]{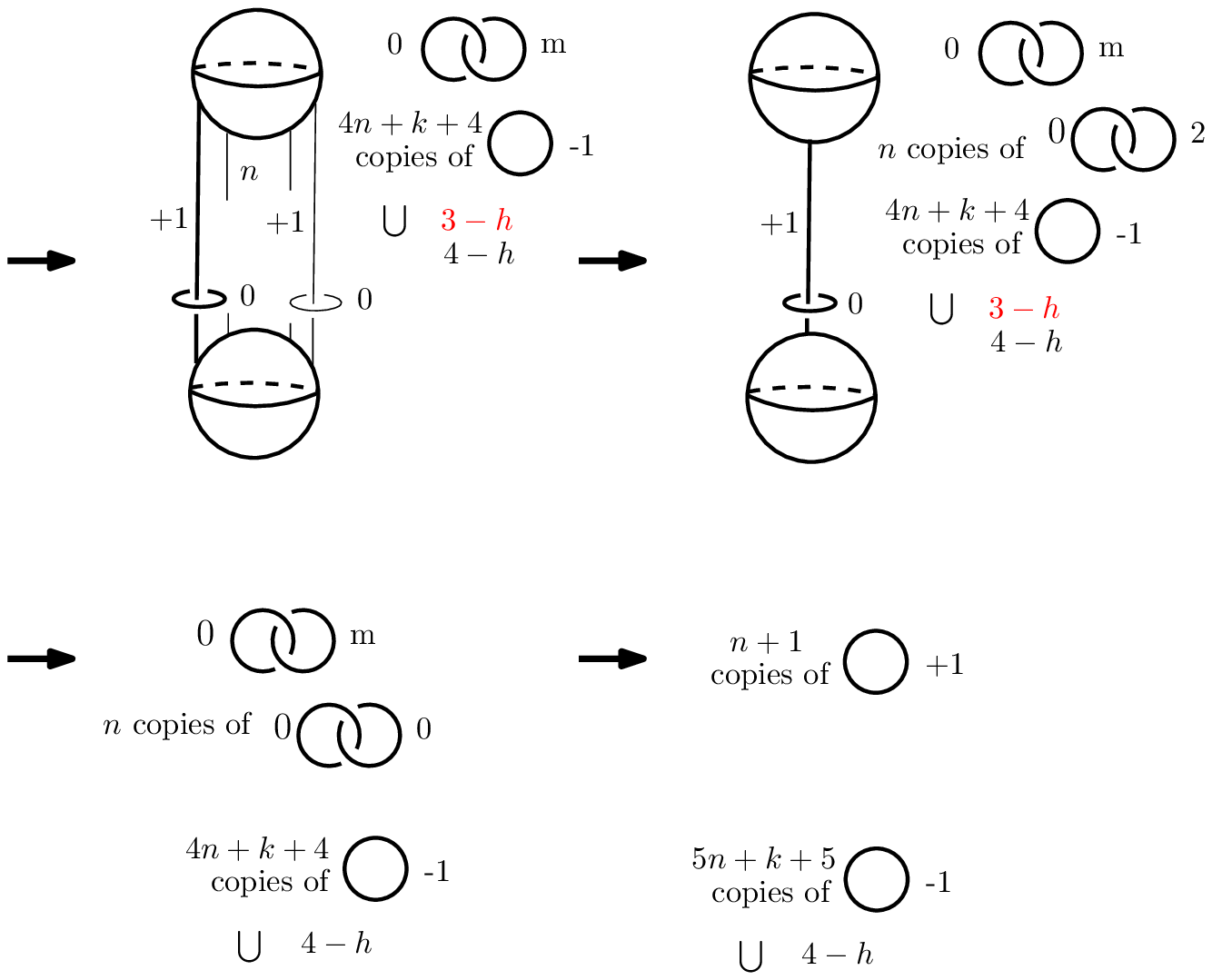}
\caption{\small }
\label{calculus2b}
\end{center}
\end{figure}

The rest of the calculus is captured in Figures \ref{calculus2a} and \ref{calculus2b}. Here any $2$-handle whose framing is not indicated should be understood to have framing $-1$. In the following paragraphs we spell out the details of this calculus, step by step. 

The first diagram in Figure \ref{calculus2a} is obtained after canceling the $1$-handle against the $0$-framed $2$-handle of the round $2$-handle going over it once, and then sliding-off the remaining lowest Lefschetz handle linking with the far left Lefschetz $2$-handle only. The framing of the latter $2$-handle now becomes $0$. In this process, we slid-off all the $k$ Lefschetz handles corresponding to $(t_a)^k$ as well, so we get a total of $k+1$ $2$-handles with framing $-1$ separated from the rest of the diagram. 

To obtain the second diagram, we slide the second Lefschetz handle from the left of the first diagram over the  $0$-framed $2$-handle on the far left. Observe that these two handles were linking with the rest of the handles in the exact same way, so the handle we slide gets separated from the bigger chunk of our diagram. Then the Lefschetz $2$-handle in the very bottom of the diagram links only with the far left $2$-handle and it can be slided-off it, turning this handle's framing to $+1$.

The third diagram contains one of the crucial steps in our inductive argument that is to follow. Namely, we slide the second (originally the third) Lefschetz \linebreak $2$-handle from the left in the second diagram over the very first $2$-handle, resulting in a $0$-framed $2$-handle that no longer links with the $1$-handle and links only with this very first $2$-handle. 

The fourth diagram is obtained by sliding all the $2$-handles that link with the far left $2$-handle, using the $0$-framed $2$-handle linking with it. In particular this splits off another $(-1)$-framed $2$-handle. Importantly, the diagram now reduces to the far left $2$-handle going over the $1$-handle once and with a $0$-framed $2$-handle that appears as a meridian to it, and $3n$ pairs of Lefschetz $2$-handles following the very same pattern that $3n+3$ pairs did so before.

The first diagram of Figure \ref{calculus2b} is simply obtained by inductively applying the above calculus for all these triples of pairs of Lefschetz $2$-handles. We therefore get a total of $4n+4$ copies of $(-1)$-framed $2$-handles disjoint from the rest of the diagram, in addition to the $k$ copies of $(-1)$-framed $2$-handles we have separated at the very beginning. 

We pass to the second diagram by sliding the $m$ copies of $(+1)$-framed $2$-handles over the far left $(+1)$-framed $2$-handle and then unlinking them using the $0$-framed $2$-handle. Doing this for all, we produce $n$ pairs of $0$ and $2$-framed $2$-handles linking once. 

The third diagram is obtained by canceling the $1$-handle against the \linebreak $(+1)$-framed $2$-handle running over it once, and then canceling the $0$-framed unlinked unknot against the $3$-handle. Moreover, the pairs of $0$ and $2$-framed $2$-handles can be turned into pairs of $0$ framed $2$-handles linking once with each other using the standard sliding argument for each one of them. (That is, slide the $2$-framed handle over the linking $0$-framed handle to change the framing of it by two, while still keeping them linked with each other only once.) 

For the fourth and last diagram, we observe that in the presence of $(-1)$-framed $2$-handles (which we will have even when $m$ and $k$ are zero), the pair of $0$ and $m$-framed $2$-handles linking once can be turned into a disjoint pair of $+1$ and $(-1)$ framed $2$-handles. Hence, the $4$-manifold which is the total space of this SBLF is 
\[(n+1) \CP \, \# \, (5n+k+5) \CPb .\] 

We can now turn to the remaining case; when the boundary fibration of the higher side is identified with the higher end of the round cobordism using a twisted gluing. When the fibration is cooked up this way, instead of the $2$-handle representing the section $S$, we get a $2$-handle linking with the $1$-handles and the \linebreak $0$-framed $2$-handle corresponding to the fiber (and linking with no other handle). However, this handle could be slid-off from the first $1$-handle (for $a$) using the \linebreak $0$-framed $2$-handle of the round $2$-handle first, and from the second $1$-handle (for $b$) using the $(+1)$-framed $2$-handle in the second diagram in Figure \ref{calculus2b}. Moreover, it can be unlinked from this $(+1)$-framed $2$-handle using the $0$-meridian repeatedly. It is not hard to see that we will end up getting an unknot that links with the $0$-framed $2$-handle (originally corresponding to the fiber) at the end. For this to prescribe a $4$-manifold they should link only once, which ends up giving us the same handlebody picture, i.e. the final diagram in Figure \ref{calculus2b}.

Lastly, let us note that $S^2 \x S^2 \# S^1 \x S^3 \# \CPb = \CP \# \CPb \# S^1 \x S^3 \# \CPb$, and $L_n \# \CPb = L'_n \# \CPb$, which can be easily seen by sliding the $2$-handles with framings $m$ and $l$ in Figures \ref{trivialgenusoneSBLFs} and \ref{trivialgenusoneSBLFs2}, respectively, over the additional $(-1)$-framed unknotted $2$-handle, and then separating this $(-1)$-framed $2$-handle off again using the $0$-framed meridians. 

We can now summarize our results. Letting the ``standard'' broken Lefschetz fibrations on respective $4$-manifolds refer to those we have considered above, we have:

\begin{thm} \label{classifymanifolds}
If $f: X \to S^2$ is a genus one relatively minimal simplified broken Lefschetz fibration, then one of the following holds:
\begin{itemize}
\item If the singular set is empty, then $(X,f)$ is $T^2 \x S^2$, $S^1 \x S^3$ or $S^1 \x L(n,1)$ equipped with standard locally trivial torus fibrations. \
\item If round singular set is not empty, but the set of Lefschetz critical points is, then $(X,f)$ is $S^2 \x S^2 \# S^1 \x S^3$,  $\CP \# \CPb \# S^1 \x S^3$, $S^4$, $L_n$, or $L'_n$, $n>1$ equipped with standard broken fibrations, respectively. \
\item If round singular set is empty, but the set of Lefschetz critical points is not, then $(X,f)$ is Hurwitz equivalent to $E(n)$ with monodromy factorization 
$(t_a, t_b)^{6n}$,    
$n>0$. \
\item If neither the round singular set nor the set of Lefschetz critical points is empty, then, either $(X,f)$ is any one of $S^2 \x S^2 \# S^1 \x S^3 \# k \CPb = \CP \# \CPb \# S^1 \x S^3 \# k \CPb$, $\# k \CPb$, $L_n \# k \CPb = L'_n \# k \CPb$, $n>1$, with higher side monodromy factorization $(t_a)^k$, or possibly after blow-ups, $(X,f)$ is $(n+1) \CP \, \# \, (5n+k+5) \CPb$ with higher side monodromy factorization Hurwitz equivalent to $(t_a, t_b)^{3n+3} \cdot (t_a)^k$, for $n, k$ non-negative integers. \
\end{itemize}
\end{thm}

Next corollary concerns solely with the topology of $4$-manifolds admitting genus one simplified broken Lefschetz fibrations:

\begin{cor} \label{classifymanifolds2}
Let $X$ be a $4$-manifold admitting a genus one (possibly not relatively minimal) simplified broken Lefschetz fibration, and $k$ be its euler characteristics. Then one of the following holds:
\begin{itemize}
\item $\pi_1(X)=\Z$ and $X=S^1 \x S^3\# k \CPb$,  $\CP \# S^1 \x S^3 \# (k-1) \CPb$ for $k>1$, or $S^2 \x S^2 \# S^1 \x S^3$.  \
\item $\pi_1(X)=\Z_n$ and $X=L_n \# k \CPb$ or $L'_n \# k \CPb$, $n>1$. \
\item $\pi_1(X)=\Z \x \Z$ and  $X=T^2 \x S^2 \# k \CPb$. \
\item $\pi_1(X)=\Z \x \Z_n$ and $X=S^1 \x L(n,1) \#k \CPb$, $n>1$.  
\item $\pi_1(X)=1$, and either $X=E(n)$ with $k=12n$, or otherwise, possibly after some blow-ups, $X$ and $n \CP \, \# \, 5n \, \CPb$, for some non-negative $n$, become diffeomorphic to each other. \
\end{itemize}
\end{cor}

\begin{proof} 
Fundamental groups and euler characteristics of the $4$-manifolds mentioned in the statement are well-known, and they can be easily computed using the handle diagrams we have given above as well. We obtain the above list after regrouping all the $4$-manifolds given in Theorem \ref{classifymanifolds}. (Here $n \CP \# 5n \, \CPb = S^4$ for $n=0$.)
\end{proof}

\vspace{0.2in}
\begin{rk} Except for $X=S^4, S^1 \x S^3, L_n, L'_n$, or their blow-ups, in all the cases covered in Theorem \ref{classifymanifolds}, the total space $X$ has $b^+(X)>0$, and therefore admits a near-symplectic form. For each case, we have depicted above a fibration with a section. Since the fibers are connected, it follows that fibers are homologically essential, and in turn, these are near-symplectic broken Lefschetz fibrations in the sense of \cite{ADK, B1, B2}.
\end{rk}

\begin{rk}
In the last case given in Theorem \ref{classifymanifolds}, we classified the ambient \linebreak $4$-manifolds only up to blow-ups. This is due to the essential role that the classification of Lefschetz fibrations over the $2$-disk with $\tilde{ \mu} = (s_1 s_2)^{3 q} s_1^{k}$ played in our proof. There we assume that $k$ is a non-negative integer,  which we achieved by passing to a blow-up of the original fibration. We shall note that these blow-ups can be performed symplectically when we have a near-symplectic broken Lefschetz fibration in hand. 

Nevertheless, there are genus one SBLFs which are not listed in our Theorem, unless one passes to a blow-up of them. A nice example, due to Tim Perutz (twisted case of which was included later in \cite{B1}), is given by the following diagram:
\begin{figure}[htp]
\begin{center}
\includegraphics[scale=1.2]{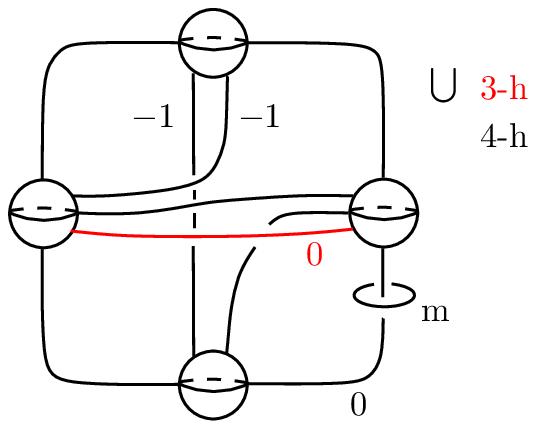}
\caption{\small Genus one SBLFs on $S^2 \x S^2$ and $\CP \# \CPb$, for $m:\text{even}$ and $\text{odd}$, respectively.}
\label{timsexample}
\end{center}
\end{figure}

\noindent The reader can verify that the global monodromy of the higher side of this fibration maps to 
$$
\tilde{\mu} = 
\left( \begin{array}{cc}
-1 & \, 0 \\
\, 4 & -1 
\end{array} \right),
$$
in $SL(2, \Z)$, which by our treatment, requires blow-ups. As shown in \cite{B1}, the total space $X$ of this fibration is $S^2 \x S^2$ or $\CP \# \CPb$, depending on the parity of $m$. (There appears to be a typo in \cite{B1}; the lower left entry of $\tilde{\mu}$ should read $4$ not $2$.) To apply our algorithm, we blow-up $X$ four times, and pass to a genus one SBLF, with  
monodromy factorization  $(s_1, s_2)^3$ of $-\text{id}$ in $SL(2, \Z)$.  
Our proof of Theorem \ref{classifymanifolds} verifies that the total space is $\CP \# 5 \CPb$, as expected.
\end{rk}

\vspace{0.1in}
\noindent \textit{Acknowledgments.} After making our preprint publicly available, Kenta Hayano informed us that he independently classified genus one simplified broken Lefschetz fibrations with at most five Lefschetz critical points in a recent preprint \cite{H}. We would like to thank Kenta Hayano for pointing out the untreated case in the proof of Theorem \ref{classifymanifolds}, Mehmetcik Pamuk and Nathan Sunukjian for their comments on a draft of this paper. The first author was partially supported by the NSF grant DMS-0906912. The second author was partially supported by JSPS KAKENHI 21340015.

\vspace{0.4in}

\end{document}